\documentclass[journal,letterpaper,twocolumn,twoside]{IEEEtran}

\usepackage{times}
\usepackage{amsfonts,amsmath,amssymb,mathrsfs,float,latexsym,amstext}
\usepackage{amsthm}
\usepackage{graphics,graphicx}
\usepackage{subfigure}
\RequirePackage[OT1]{fontenc}
\usepackage{cite}

\usepackage[colorlinks=true,linkcolor=black,anchorcolor=black,citecolor=black,filecolor=black,menucolor=black,runcolor=black,urlcolor=black]{hyperref}

\newcommand{\ignore}[1]{} 

\allowdisplaybreaks 

\theoremstyle{plain}
\newtheorem{theorem}{Theorem}

\newtheorem{lemma}{Lemma}

\theoremstyle{definition}

\newtheorem{remark}{Remark}

\font\msbmx=msbm10                   
\font\msbmvii=msbm7                  
\font\msbmv=msbm5
\newcommand{\mb}[1]{\mathbf{#1}} 

\newcommand{\Xb}{{\mb{X}}}

\newcommand{\Yb}{{\mb{Y}}}

\newcommand{\mbs}[1]{\boldsymbol{#1}} 

\newcommand{\mc}[1]{\mathcal{#1}} 

\newcommand{\Nc}{{\mc{N}}}

\newcommand{\Ac}{{\mc{A}}}

\newcommand{\Ic}{{\mc{I}}}

\usepackage{mathrsfs}
\newcommand{\Pb}{{\mathsf{P}}} 

\newcommand{\Eb}{{\mathsf{E}}}

\newcommand{\LPFA}{\mathsf{LPFA}}
\newcommand{\Hyp}{{\mathsf{H}}} 
\newcommand{\PD}{\mathsf{PD}}

\newcommand{\mrm}[1]{\mathrm{#1}}

\newcommand{\mbb}[1]{\mathbb{#1}} 
\def\One{\mathchoice{\rm 1\mskip-4.2mu l}{\rm 1\mskip-4.2mu l}
{\rm 1\mskip-4.6mu l}{\rm 1\mskip-5.2mu l}}
\newcommand\Ind[1]{{\One_{\{#1\}}}} 

\newcommand{\Zbb}{\mbb{Z}} 

\newcommand{\class}{{\mbb{C}}}

\newcommand{\Tcs}{T_{\mbox{\rm\tiny CS}}}

\newcommand{\Tfma}{T_{\mbox{\rm\tiny FMA}}}

\newcommand{\set}[1]{\left\{#1\right\}}

\newcommand{\brc}[1]{\left(#1\right)}
\newcommand{\brcs}[1]{\left[#1\right]}

\renewcommand{\le}{\leqslant} 
\renewcommand{\ge}{\geqslant}
\newcommand{\esssup}{\operatornamewithlimits{ess\,sup}}
\newcommand{\essinf}{\operatornamewithlimits{ess\,inf}}
\def\argmin{\mathop{\rm arg\,min}\limits}

\usepackage{color}
\definecolor{OliveGreen}{rgb}{.2,0.6,0.2}
\definecolor{BrickRed}{rgb}{.7,0.2,0.2}
\definecolor{Violet}{rgb}{.58,0.09,0.65}

\begin{document}

\title{Optimal Sequential Detection of Signals with Unknown Appearance and Disappearance Points in Time
\thanks{The work was supported in part by the grant 18-19-00452
from the Russian Science Foundation at the Moscow Institute of Physics and Technology.}}

\author{Alexander~G.~ Tartakovsky, \IEEEmembership{Senior~Member,~IEEE},
\thanks{A. G. Tartakovsky is a Deputy Head of the Space informatics Laboratory at the Moscow Institute of Physics and Technology, Russia
 and President of AGT StatConsult, Los Angeles, California, USA; e-mail: agt@phystech.edu} Nikita R. Berenkov,
 \thanks{N. R. Berenkov is a postgraduate student at the Moscow Institute of Physics and Technology, Russia; email: nberenkov@mail.ru}
Alexei E. Kolessa, \thanks{A. E. Kolessa is the Principal Scientist in the Space informatics Laboratory at the Moscow Institute of Physics and Technology, Russia;
e-mail: kolessa.ae@phystech.edu}  and Igor~V.~Nikiforov \thanks{I. V. Nikiforov is Emeritus Professor at the
Universit\'e de Technologie de Troyes, Troyes Cedex, France; e-mail: Igor.Nikiforov@utt.fr}
\thanks{Manuscript received 2020; revised 2021} 
\thanks{Copyright (c) 2021 IEEE. Personal use of this material is permitted.  However, permission to use this material for any other purposes must be obtained from the
IEEE by sending a request to pubs-permissions@ieee.org.}
}

\markboth{IEEE Transactions on Signal Processing,~Vol.~~, No.~~, ~~~2021}%
{Tartakovsky, Berenkov, Kolessa, Nikiforov: Optimal Sequential Detection of Signals}

\maketitle

\begin{abstract}
The paper addresses a sequential changepoint detection problem, assuming that the duration of change may be finite and unknown. This problem is of importance for many
applications, e.g., for signal and image processing where signals appear and disappear at unknown points in time or space. In contrast to the conventional optimality criterion in quickest change detection that requires minimization of the expected delay to detection for a given average run length to a false alarm, we focus on a reliable maximin change detection criterion of maximizing the minimal probability of detection in a given time (or space) window for a given local maximal probability of false alarm in the prescribed window. We show that the optimal detection procedure is a modified CUSUM procedure. We then compare operating characteristics of this optimal procedure with popular in engineering the Finite Moving Average (FMA) detection algorithm and the ordinary CUSUM procedure using Monte Carlo simulations, which show that typically the later algorithms have almost the same performance as the optimal one. At the same time, the FMA procedure has a substantial advantage -- independence to the intensity of the signal, which is usually unknown. Finally, the FMA algorithm is applied to detecting faint streaks of satellites in optical images.
\end{abstract}




%
\section{Introduction} \label{sec:intro}

\IEEEPARstart{C}hangepoint detection problems arise in a variety of applications that are described in detail in \cite{TNB_book2014}. In most cases, the quickest change detection is considered where one has to
detect a change as soon as possible, i.e., with the minimal average delay to detection for a given false alarm rate (see, e.g.,\cite{TNB_book2014} and references therein). However, in certain applications, it is of interest to consider
a reliable change detection when it is needed to maximize a probability of detection at a certain time (or space) interval for a given probability of false alarm.
For instance, in surveillance systems such as radars, sonars, and electro-optic/infrared sensor systems, which deal with detecting moving and maneuvering objects that appear and disappear at unknown times, it is necessary to detect
a signal from a randomly appearing target in clutter and noise with the maximal detection probability \cite{Bakutetal-book63,Richards_radar2014,Marage_sonar2013,Tartakovsky&Brown-IEEEAES08}. Also, examples include but are not
limited to: (a) aircraft navigation with many safety-critical modes (landing, takeoff, etc.) \cite{BakNik00}, where the minimum operational performance specifies the required time-to-alert, the worst-case missed detection probability and
the worst-case probability of false alarm during a given period; and (b) cyber-security \cite{Do2017,Tartakovsky-Cybersecurity14,Tartakovskyetal-SM06,Tartakovskyetal_IEEESP2013} when there are malicious intrusion attempts in computer networks which incur significant financial damage and cause severe harm to the integrity of personal information. In these cases, it is essential to devise effective techniques to detect anomalies in observations reliably so that an
appropriate response can be provided and the negative consequences are mitigated. In these and other applications, the statistical behavior of observed data is dynamic, so it is of importance to detect transient changes. For example,
after an outage in the power systems, the system's transient behavior is dominated by the generators' inertial response.

In this paper, we address the problem of detecting a change that has a finite duration. The problem of detecting transient
changes with known and unknown durations has been considered in~\cite{Egea-Roca-tsp2018,NikiforovetalSQA2012, NikiforovetalIEEEIT2017, NoonZhig20, Repin-pit91,Tartakovsky-pit87,Tartakovsky-pit88a}. In particular, articles
\cite{NikiforovetalSQA2012, NikiforovetalIEEEIT2017} establish the asymptotic performance of the window-limited CUSUM procedure as the false alarm probability goes to $0$
for detection of transient changes with known
duration. However, the issue of optimality or asymptotic optimality is still open. The problem of detection of transient and moving anomalies has also been considered in papers
\cite{Premkumar2010BayesianQT, Tchamkerten2017, Rovastos_2017, Zou_2017, Rovatsos-sa20} but in terms of quickest change detection.

The rest of the paper is organized as follows. In Section~\ref{sec:Model}, we describe the stochastic model, which is treated in the paper, as well as the optimality criteria.
In Section~\ref{sec:Rule}, we find the optimal detection procedure that maximizes detection probability in the worst-case scenario in the class of detection procedures with the given local false alarm
probability (in a certain window), assuming that the duration of a change (or window size) is random and distributed with the geometric distribution. This procedure turns out to be the modified 
Cumulative Sum (CUSUM) rule. 
Proof of optimality (see Theorem~\ref{Th:MaximinOptimalityMm}) is based on the optimal stopping results established in the paper by Poor \cite{poor-as98}.
In Section~\ref{sec:ADP}, two alternative detection procedures are introduced -- the Finite Moving Average (FMA) procedure and the conventional CUSUM procedure.
In Section~\ref{sec:MC}, we provide the results of Monte Carlo simulations for the Gaussian model which show that the FMA and the CUSUM procedures have almost the same operating
characteristics as the optimal procedure. This allows us to suggest using the FMA procedure in Section~\ref{sec:Appl} in an important practical problem of detecting streaks of space objects with
unknown position (beginning and end) in 2-D images obtained by telescopes. Section~\ref{sec:Remarks} concludes the paper.

\section{The Stochastic Model and Optimality Criteria}\label{sec:Model}

Suppose there is a sequence of independent observations $\{Y_n\}_{n \ge 1}$, observed sequentially in time subject to a change at an unknown time $\nu\in\{0,1, 2, \dots\}$, which lasts till the time $\nu+ N$
so that $Y_1,\dots,Y_\nu$ and $Y_{N+\nu+1}, Y_{N+\nu+2}, \dots$ are generated by one stochastic model and $Y_{\nu+1}, Y_{\nu+2}, \dots, Y_{\nu+N}$ by another model.
Throughout the paper, $\nu$ is treated as unknown and nonrandom, and $N$
can be either unknown deterministic (as in Section~\ref{sec:Appl}) or random (as in Section~\ref{sec:Rule}). See also Remark~\ref{Rem:Thconstraints}.
The joint density $p(\Yb_1^n | \Hyp_{\nu,N})$ of the vector
$\Yb_1^n= (Y_1,\dots,Y_n)$ observed up to time $n$ under the hypothesis $\Hyp_{\nu,N}$ that the change happens at the time $\nu$ and ends at $N$ is of the form
\begin{equation}\label{modeldetisol}
\begin{split}
& p(\Yb_1^n | \Hyp_{\nu,N}) = p(\Yb_1^n| \Hyp_\infty)= \prod_{t=1}^n g(Y_t)
\\
& ~~\text{for}~~ \nu \ge n ,
\\
& p(\Yb_1^n | \Hyp_{\nu,N}) = \prod_{t=1}^{\nu} g(Y_t) \times \prod_{t=\nu+1}^{n} f(Y_t)
\\
& ~~\text{for}~~ \nu <n \le \nu+ N,
\\
& p(\Yb_1^n | \Hyp_{\nu,N}) = \prod_{t=1}^{\nu} g(Y_t) \times \prod_{t=\nu+1}^{\nu+N} f(Y_t)
\\
& \times \prod_{t=\nu+N+1}^n g(Y_t) ~~ \text{for}~~ n > \nu+ N ,
\end{split}
\end{equation}
where $g(Y_t)$ and $f(Y_t)$ are pre- and post-change densities, respectively. The event $\{\nu=\infty\}$ and the
corresponding hypothesis $\Hyp_\infty: \nu=\infty$ mean that there never is a change. Notice that the model \eqref{modeldetisol} implies that $Y_{\nu+1}$ is the first post-change observation
under hypothesis $\Hyp_{\nu,N}$.

A sequential changepoint detection procedure $T$ is a stopping time associated with the time of alarm on change.

Conventional quickest detection optimality criteria require minimizing the average delay to detection for a given false alarm rate at an
infinite time horizon (assuming $N=\infty$) and do not consider a probability of detection of a change in a given fixed time interval \cite{ShiryaevTPA63, ShiryaevBook78, lorden-ams71,PollakAS85,TNB_book2014}.
Often, however, practitioners are interested in such probabilities under a given false alarm rate even if the change lasts infinitely
long.\footnote{In practice, this means that the length of a change is substantially larger than an average detection delay.} Besides, in many applications, the length of a change $N$ is finite, e.g., in
problems of detecting transient
changes with known and unknown durations~\cite{NikiforovetalSQA2012, NikiforovetalIEEEIT2017, Repin-pit91,Tartakovsky-pit87,Tartakovsky-pit88a}.
Then stopping outside of the interval $(\nu, \nu+ N]$ of a given duration $N$ may not be quite appropriate. For example, in the context of safety-critical systems, serious degradation of the system
security occurs when the transient change is detected with a delay greater than a required time-to-alert. Therefore, the probability of detection of a change within a given fixed time interval should be used instead
of the average delay to detection. In such cases, it is reasonable to find detection rules that maximize the probability of
detection in a certain fixed time interval $(\nu, \nu +M], M \le N$, and to consider the following optimality criterion:
{\em Find a rule $T_{\mathrm{opt}}\in\class(m,\alpha)$ such that for every $0<\alpha < 1$ and some $m \ge 1$}
\begin{equation} \label{MaximinPDOpt1Mrandom}
\begin{split}
& \inf_{\nu\in \Zbb_+} \essinf\, \Pb_\nu(T_{{\rm opt}} \le \nu+M| \Yb_1^\nu, T_{{\rm opt}} >\nu)
 \\
 &=\sup_{T\in\class(m, \alpha)} \inf_{\nu\in \Zbb_+} \essinf\, \Pb_\nu(T \le \nu+M | \Yb_1^\nu, T >\nu) ,
 \end{split}
\end{equation}
where
\[
\class(m, \alpha)= \set{T: \sup_{\ell \in \Zbb_+} \Pb_\infty(T \le \ell+m | T> \ell) \le \alpha}
\]
is the class of detection procedures for which the local maximal probability of false alarm
\[
\LPFA_m(T)=\sup_{\ell \in \Zbb_+} \Pb_\infty(T \le \ell+m | T> \ell)
\]
 in a time interval of a fixed length $m \ge 1$ does not exceed a predefined level $\alpha\in(0,1)$.\footnote{In general, $m$ and $M$ are different.} 
 Hereafter $\Pb_\nu$ denotes the probability under which the change occurs at $\nu \in \Zbb_+$ and $\Pb_\infty$ when the change never happens;  $\Zbb_+=\{0,1,2,\dots\}$ states for the set of 
 nonnegative integers. Also, $\essinf$ and $\esssup$ denote essential infimum and essential supremum, respectively.

Solving the optimization problem \eqref{MaximinPDOpt1Mrandom} for any fixed $M$ (except for $M=1$) is very difficult (see \cite{Moustakides14,PollakKriegerSQA2013,Tartakovsky_book2020}
for some optimal properties established for the case $M=1$). Even in an asymptotic setting as $\alpha\to0$ this problem is still open.
Assume now that $M$, $M \le N$ is finite and random with the given distribution $\pi_i=\Pr(M=i)$, $i =1,2,\dots$. In particular, this assumption is reasonable when $M=N$ and the unknown duration of a change $N$ is a nuisance parameter, i.e., the fact of change disappearance does not have to be detected. The results also valid when the change is persistent, i.e., $N=\infty$, 
but still, the goal is to maximize the probability of detection in a time window of a random size $M$.

Introduce the probability measure
$$
\bar{\Pb}_\nu(\Ic\times\Ac) = \sum_{i\in \Ic} \pi_i \Pb_\nu(\Ac | M=i).
$$

Then the probability of detection is re-written as
\begin{equation}\label{PDaverage}
\begin{split}
& \bar{\Pb}_\nu(T\le \nu+M|\Yb_1^\nu,T>\nu)
\\
&= \sum_{i=1}^\infty \pi_i \Pb_\nu(T\le \nu+i|\Yb_1^\nu,T>\nu, M=i),
\end{split}
\end{equation}
and the optimality criterion \eqref{MaximinPDOpt1Mrandom} gets modified as
\begin{equation} \label{MaximinPDOpt1Mrandom_1}
\begin{split}
 &\inf_{\nu\in \Zbb_+} \essinf\, \bar{\Pb}_\nu(T_{{\rm opt}} \le \nu+M| \Yb_1^\nu, T_{{\rm opt}} >\nu)
 \\
 &=\sup_{T\in\class(m, \alpha)} \inf_{\nu\in \Zbb_+} \essinf\, \bar{\Pb}_\nu(T \le \nu+M | \Yb_1^\nu, T >\nu).
 \end{split}
\end{equation}

In the next section, we provide a solution to this problem for the geometric distribution $\pi_i$.

\section{An Optimal Detection Procedure}\label{sec:Rule}

Let $\pi_i$ be the geometric distribution $\mrm{Geom}(\varrho)$ with the parameter $\varrho\in(0,1)$:
\[
\pi_i =\varrho(1-\varrho)^{i-1}, \quad i =1,2,\dots.
\]
Let $\Lambda_n=f(Y_n)/g(Y_n)$ ($n=1,2,\dots$) be the likelihood ratio and introduce the statistic $V_\varrho(n)$ by the recursion
\begin{equation}\label{Vrho0}
V_{\varrho}(n) = \max\set{1,V_{\varrho}(n-1)} \Lambda_{n} (1-\varrho), \quad n \ge 1
\end{equation}
with the initial condition $V_{\varrho}(0) = 1$ as well as the associated stopping time
\begin{equation}\label{MaximinOptruleMm}
T_{\varrho}(B)=\inf\set{n \ge 1: V_{\varrho}(n) \ge B}, \quad B >0.
\end{equation}
This kind of statistic first appears in the paper by Poor \cite{poor-as98} who considered the exponential delay function in the minimax ``quickest'' detection problem. Also, his results will be used in the 
proof of Theorem~\ref{Th:MaximinOptimalityMm}.

Let $\Eb_{\nu}$ and $\Eb_\infty$ denote expectations under probability measures $\Pb_\nu$ and $\Pb_\infty$, respectively, where $\Pb_{\nu}$ corresponds to
model \eqref{modeldetisol} with an assumed value of the change point $\nu$. 

The following theorem, whose proof is given in the Appendix, establishes the structure of the optimal detection procedure.

\begin{theorem}\label{Th:MaximinOptimalityMm}
Let observations $\{Y_n\}_{n\ge 1}$ be independent with a density $g(x)$ pre-change and with a density $f(x)$ post-change. Suppose the distribution of the window size $M$ is
$\mrm{Geom}(\varrho)$. Further, assume that the $\Pb_\infty$-distribution of the likelihood ratio $\Lambda_1=f(Y_1)/g(Y_1)$ is continuous and that $\Pb_\infty\{\Lambda_1> (1-\varrho)^{-1}\} =1$. Then the change detection rule $T_{\varrho}(B)$ defined in
\eqref{MaximinOptruleMm} with the statistic $V_{\varrho}(n)$ given by the recursion \eqref{Vrho0}
and with threshold $B=B_{m,\alpha}$ that satisfies
\begin{equation}\label{Bmalpha}
\sup_{\ell \in \Zbb_+} \Pb_\infty \set{T_{\varrho}(B) \le \ell +m | T_{\varrho} > \ell} = \alpha
\end{equation}
is maximin optimal in the problem \eqref{MaximinPDOpt1Mrandom_1} for all $0<\alpha <1$.
\end{theorem}

Note that statistic $V_\varrho(n)$ is the maximal weighted likelihood ratio
\[
V_\varrho (n) =\max_{1 \le k \le n} \brcs{ \prod_{j=k}^n (1-\varrho) \Lambda_j},
\]
so the optimal rule \eqref{MaximinOptruleMm} is nothing but a modified CUSUM rule with an additional factor $1-\varrho$.\footnote{In the standard CUSUM procedure $\varrho =0$.}
If the distribution $\Pb_\infty(\Lambda_1 \le y)$ is not continuous the assertion of Theorem~\ref{Th:MaximinOptimalityMm} holds for a randomized procedure with a randomization on the boundary $B$.

\begin{remark}
As follows from the proof, the detection algorithm \eqref{Vrho0}--\eqref{MaximinOptruleMm} is also optimal in the class of procedures subject to the constraint on the Average Run Length to False Alarm (ARL2FA) $\Eb_\infty[T] \ge \gamma$ since class $\class(m,\alpha)$ is more stringent than $\class_\gamma=\{T: \Eb_\infty[T]\ge \gamma\}$ for some appropriately selected $\gamma=\gamma(m,\alpha)$ (see Lemma~\ref{Lem:LCPFAtoARL} in the Appendix). However, ARL2FA makes sense only if the $\Pb_\infty$-distribution of stopping times of detection procedures is geometric or close to geometric -- asymptotically exponential. Asymptotic exponentiality property holds for many detection procedures with Markov detection statistics \cite{PollakTartakovskyTPA09}. However, we do not know whether this is correct for the FMA procedure considered below. If the no-change distribution of stopping times is not close to geometric, then a large value of ARL2FA does not guarantee a small value of the maximal local false alarm probability $\LPFA_m(T)$, which is usually a necessary property in practice. A detailed discussion of this issue may be found in \cite{Tartakovsky-SQA08a}.
\end{remark}

\begin{remark} \label{Rem:Thconstraints}
The assertions of Theorem~\ref{Th:MaximinOptimalityMm} hold in two cases: (a) when the window size is equal to the unknown change duration, $M=N\sim \mrm{Geom}(\varrho)$, and (b) when the change is persistent, $N=\infty$, and $M\sim \mrm{Geom}(\varrho)$. However, the latter case has perhaps only theoretical rather than practical significance.
\end{remark}

\begin{remark}
The time index $n$ in all previous formulas can be replaced by the argument of the vector-valued function of position $(x_i,y_j)$,
as it is done in Section~\ref{sec:Appl} for the problem of detecting objects in two-dimensional images.
\end{remark}

\section{Alternative Detection Procedures} \label{sec:ADP}

In this and the next section, we set $M=N$, i.e., we assume that the window size $M$ is equal to the signal duration $N$.

While detection procedure $T_{\varrho}$ given by \eqref{Vrho0}--\eqref{MaximinOptruleMm} is strictly optimal, it requires a strong assumption on the geometric distribution of the signal duration $N$.
A more practical approach is to use the procedures in a fixed sliding window with size $L$. In papers \cite{NikiforovetalSQA2012, NikiforovetalIEEEIT2017},
a window-limited CUSUM procedure
\begin{equation}\label{WLCUSUM}
\overline{T}=\inf\left\{n\geq L : \max\limits_{1 \leq k \leq L}\left[\sum_{t=n-k+1}^{n}\lambda_t-A(k)\right]\geq 0\right\},
\end{equation}
where $\lambda_t= \log \Lambda_t$ is the log-likelihood ratio, has been proposed for this purpose. It was shown that
the Finite Moving Average procedure given by the stopping time
\begin{equation}\label{FMA}
\Tfma(a) = \inf \left\{n \ge L: \sum_{t = n - L + 1}^{n} \lambda_t \ge a\right\}
\end{equation}
and the window-limited CUSUM procedure $\overline{T}$ with a specific (optimal) variable threshold $A(k)$, minimizing
the worst-case missed detection probability, have the same asymptotic  performance as the maximal probability of false alarm $\alpha\to0$.
Window-limited procedures were also considered by Lai in \cite{LaiIEEE98} and shown to be asymptotically optimal in the quickest change detection problem with persistent changes for minimizing average detection delay for i.i.d. and non-i.i.d stochastic models.

Another reasonable method is a simple CUSUM procedure. It is easy to show that maximizing the likelihood ratio over the unknown points of change appearance and disappearance  ($\nu$ and $N$, respectively) leads to Page's CUSUM statistic:
\begin{equation}
\begin{split}
&V(n) = \max_{\nu \ge 0}\max_{N \ge 1} \frac{p(\Yb_1^n | \Hyp_{\nu,N})}{p(\Yb_1^n| \Hyp_\infty)} \\
& = \max\set{1,V(n-1)} \Lambda_{n}, \quad n \ge 1.
\end{split}
\end{equation}
Hence, define the CUSUM procedure by the stopping time:
\[
\Tcs(C) =\inf\set{n \ge 1: V(n) \ge C}, \quad C >0.
\]

\begin{remark} \label{Rem:Monotone}
If the LLR $\lambda_t$ is a monotone function of the statistic $S_t$, then the FMA procedure can be written as
\begin{equation}\label{TFMA}
\Tfma(\tilde{a}) = \inf \left\{n \ge L: \sum_{t = n - L + 1}^{n} S_t \ge \tilde{a}\right\}.
\end{equation}
If the post-change distribution depends on an unknown parameter $\theta$, then $\lambda_t=\lambda_t(\theta)$ and the optimal modified CUSUM $T_{\varrho}(B)$ as well as the ordinary CUSUM $\Tcs(C)$
depend on $\theta$. Therefore, they are sensitive to the mismatch of the true value and  assumed values of $\theta$, while the FMA structure does not depend on $\theta$.
We expect that the FMA procedure maximizes  (approximately)
the probability of detection $\PD_\theta$ uniformly for all parameter values.
\end{remark}

In the next section, we compare the performance of detection procedures $\Tfma(\tilde{a})$ and $\Tcs(C)$  with the optimal one.

\section{Monte Carlo Simulations}\label{sec:MC}

We stress that in simulations the time window size $M=N$ is assumed random with the geometric distribution, $N \sim \mrm{Geom}(\varrho)$. The window's length $L$ in the FMA
rule is fixed and selected as $L = \Eb[N] = 1/\varrho$.

Consider the standard signal-plus-noise model
\[
Y_n = \theta \Ind{\nu < n \le \nu + N} + \xi_n, \quad n \ge 1,
\]
where $\Ind{\Ac}$ is the indicator of the event $\Ac$ and $\{\xi_{n}\}_{n\ge 1}$ is the i.i.d.\ Gaussian sequence with mean zero and standard deviation $\sigma>0$, $\xi_n \sim \Nc(0,\sigma^2)$, $\sigma = 1$.
Thus, the observations $Y_n$ have normal distribution $\Nc(0,1)$ pre-change and normal $\Nc(\theta,1)$ post-change.

The following notation is used for the minimal probability of detection of rule $T$:
	\begin{align*}
 \PD(T)& = \inf_{\nu\in \Zbb_+} \essinf\, \bar{\Pb}_\nu(T \le \nu+N| \Yb_1^\nu, T >\nu) .
 \end{align*}

In simulations, we used the following experimentally proved conjectures for FMA $T=\Tfma$ and modified CUSUM $T=T_{\varrho}$ rules:
	\begin{align*}
	\LPFA_m(T) = \Pb_{\infty}(T \le m), \quad
	\PD(T)= \bar{\Pb}_0(T \le N),
	\end{align*}
i.e., $\ell= 0$ delivers the maximum of the local PFA and $\nu = 0$ delivers minimum to the average detection probability over the distribution of $N$.

\subsection{Comparison of FMA and Modified CUSUM Rules}

In simulations, we use the following parameters:	
	\begin{enumerate}
		\item The mean after changepoint: $\theta = 2.0, 1.2$.
		\item The tuning parameter $\varrho$ of the modified CUSUM rule: $0.2, 0.1, 0.05$ ($L = 5, 10, 20$, respectively).
		\item The local PFA $\LPFA_m(T)$: $0.001$.
		\item The window length $m$ for the local PFA: $20, 80$.
    \item The number of Monte Carlo repetitions: $5\cdot 10^5$.
	\end{enumerate}
		
The results of comparing the performance of the optimal modified CUSUM and the FMA rules are shown in Table~\ref{t:B_gauss_inf_1} in the case where both rules are tuned to the same true 
parameter value $\theta$.
	
	\begin{table}[!ht!]
		\renewcommand{\arraystretch}{1}
		\begin{center}
			\caption{Detection characteristics of the modified CUSUM rule and the FMA rule\label{t:B_gauss_inf_1}}
			
			\begin{tabular}{ |c|c|c|c| }
				\hline
				\multicolumn{4}{|c|}{$\theta = 2.0, m = 20$ } \\
				\hline
				$\varrho $ & 0.2 & 0.1 & 0.05 \\
				\hline
				$\PD(T_{\varrho})$ & 0.3677 & 0.6099  & 0.7843 \\
				\hline
				$\PD(\Tfma)$ & 0.3512 & 0.5014 & 0.6394 \\
				\hline
				\hline
				\multicolumn{4}{|c|}{$\theta = 2.0, m = 80$ } \\
				\hline
				\hline
				$\PD(T_{\varrho})$ & 0.3290  & 0.5659  & 0.7797\\
				\hline
				$\PD(\Tfma)$ & 0.3220 & 0.4763 & 0.6029 \\
				\hline
				\hline
				\multicolumn{4}{|c|}{$\theta = 1.2, m = 20$ } \\
				\hline
				$\varrho $ & 0.2 & 0.1 & 0.05 \\
				\hline
				$\PD(T_{\varrho})$ & 0.1510  & 0.3547 & 0.5641 \\
				\hline
				$\PD(\Tfma)$ & 0.1424 & 0.3142 & 0.4738 \\
				\hline
				\hline
				\multicolumn{4}{|c|}{$\theta = 1.2, m = 80$ } \\
				\hline
				$\varrho$ & 0.2 & 0.1 & 0.05 \\
				\hline
				$\PD(T_{\varrho}) $ & 0.1197  & 0.2917 & 0.4910 \\
				\hline
				$\PD(\Tfma)$ & 0.1016 & 0.2694 & 0.4423 \\
				\hline
			\end{tabular}
		\end{center}
	\end{table}

It is seen from Table~\ref{t:B_gauss_inf_1} that the modified CUSUM rule outperforms the FMA rule in this setting, as expected since by Theorem~\ref{Th:MaximinOptimalityMm} it is strictly optimal.
However, the difference in performance is not dramatic. It is small for $\varrho=0.2$ and $0.1$.

For practical purposes, it is of interest to compare the performance of the modified CUSUM rule against the FMA rule under model mismatch. There are two types of mismatch: (a) the duration $N$ 
of a transient change differs from the duration defined in the detection algorithm, and (b) the true signal amplitude $\theta_r$ differs from the assumed value $\theta$ used in the detection algorithm. 

We first compare
the modified CUSUM rule against the FMA rule under the transient change duration mismatch, assuming that $N$ is fixed. The tuning parameter of the modified CUSUM rule is given as $\varrho=1/N$ 
and the FMA window's length is $L=N$. The results of this comparison are given in Table~\ref{t:B_gauss_inf_fixed}. Now, the FMA rule outperforms the modified CUSUM rule. 
However, the modified CUSUM rule retains the competitive characteristics.

\begin{table}[!ht!]
	\renewcommand{\arraystretch}{1}
	\begin{center}
		\caption{Detection characteristics of the modified CUSUM and FMA rules for mismatch in change duration \label{t:B_gauss_inf_fixed}}
		
		\begin{tabular}{ |c|c|c|c| }
			\hline
			\multicolumn{4}{|c|}{$\theta = 2.0, m = 20$ } \\
			\hline
			$N $ & 5 & 10 & 20 \\
			\hline
			$\PD(T_{\varrho})$ & 0.6739 & 0.9790  & 0.998 \\
			\hline
			$\PD(\Tfma)$ & 0.7454 & 0.9956 & 0.999 \\
			\hline
			\multicolumn{4}{|c|}{$\theta = 2.0, m = 80$ } \\
			\hline
			\hline
			$N $ & 5 & 10 & 20 \\
			\hline
			$\PD(T_{\varrho})$ & 0.5574  & 0.9632  & 0.998\\
			\hline
			$\PD(\Tfma)$ & 0.6246 & 0.9880 & 0.999 \\
			\hline
			\multicolumn{4}{|c|}{$\theta = 1.2, m = 20$ } \\
			\hline
			\hline
			$N$ & 5 & 10 & 20 \\
			\hline
			$\PD(T_{\varrho})$ & 0.0886  & 0.4662 & 0.9205 \\
			\hline
			$\PD(\Tfma)$ & 0.1355 & 0.5452 & 0.9629 \\
			\hline
			\hline
			\multicolumn{4}{|c|}{$\theta = 1.2, m = 80$ } \\
			\hline
			$N$ & 5 & 10 & 20 \\
			\hline
			$\PD(T_{\varrho}) $ & 0.0383  & 0.3356  & 0.8666 \\
			\hline
			$\PD(\Tfma)$ & 0.0739 & 0.4023 & 0.9203 \\
			\hline
		\end{tabular}
	\end{center}
\end{table}

Second, we consider the parameter $\theta$ mismatch, i.e., we compare the performance of the modified CUSUM and FMA rules when the true
signal intensity $\theta_r$ differs from the assumed value $\theta$. Specifically, the true model has the form
\[
Y_n = \theta_r \Ind{ \nu < n \le \nu + N} + \xi_n, \quad n \ge 1,
\]
where $\theta_r$ is true signal intensity. The results of this comparison are given in Table~\ref{t:mism} for the following parameters: $\varrho = 0.1,~ N \sim \mrm{Geom}(\varrho),~ L = \Eb[N] = 10,~ \LPFA_m = 0.001,~ m = 20 $ and the assumed signal intensity values are $\theta = 2.0, 1.8, 1.6, 1.4, 1.2$.

\begin{table}[!ht!]
	\renewcommand{\arraystretch}{1}
	\begin{center}
		\caption{Detection characteristics of the modified CUSUM and FMA rules with model mismatch in signal intensity \label{t:mism}}
		
		\begin{tabular}{ |c|c|c|c|c|c|c| }
			\hline
			\multicolumn{7}{|c|}{$\theta = 2.0$ } \\
			\hline
			$\theta_r$ & 1.0 & 1.2 & 1.4 & 1.6 & 1.8 & 2.0 \\
			\hline
			$\PD(T_{\varrho})$ & 0.1654  & 0.2742   & 0.3870 & 0.4772 & 0.5517 & 0.6091 \\
			\hline
			$\PD(\Tfma)$ & 0.2452 & 0.3341 & 0.3994 & 0.4436 & 0.4847 & 0.5162 \\
			\hline
			\hline
			\multicolumn{7}{|c|}{$\theta = 1.8$ } \\
			\hline
			$\theta_r$ & 1.0 & 1.2 & 1.4 & 1.6 & 1.8 & 2.0 \\
		\hline
		$\PD(T_{\varrho})$ & 0.1776  & 0.2929   & 0.3952 & 0.4852 & 0.5523 & 0.6064 \\
		\hline
		$\PD(\Tfma)$ & 0.2464 & 0.3351 & 0.3972 & 0.4465 &0.4837 & 0.5170 \\
			\hline
			\hline
		\multicolumn{7}{|c|}{$\theta = 1.6$ } \\
		\hline
		$\theta_r$ & 1.0 & 1.2 & 1.4 & 1.6 & 1.8 & 2.0 \\
		\hline
		$\PD(T_{\varrho})$ & 0.1975  & 0.3077   & 0.4088 & 0.4888 & 0.5530 & 0.6029 \\
		\hline
		$\PD(\Tfma)$ & 0.2467 & 0.3338 & 0.3992 & 0.4461 &0.4846 & 0.5178 \\
			\hline
			\hline
		\multicolumn{7}{|c|}{$\theta = 1.4$ } \\
		\hline
		$\theta_r$ & 1.0 & 1.2 & 1.4 & 1.6 & 1.8 & 2.0 \\
		\hline
		$\PD(T_{\varrho})$ & 0.2099 & 0.3193  &0.4109 & 0.4851 & 0.5454 & 0.5925 \\
		\hline
		$\PD(\Tfma)$ & 0.2453 & 0.3349 & 0.3993 & 0.4462 &0.4854 & 0.5168 \\
		\hline
		\hline
				\multicolumn{7}{|c|}{$\theta = 1.2$ } \\
			\hline
			$\theta_r$ & 1.0 & 1.2 & 1.4 & 1.6 & 1.8 & 2.0 \\
			\hline
			$\PD(T_{\varrho})$ & 0.2219 & 0.3226  &0.4083 & 0.4792 & 0.5352 & 0.5806\\
			\hline
			$\PD(\Tfma)$ & 0.2473 & 0.3357 & 0.3986 & 0.4471 &0.4847 & 0.5160 \\
		\hline
			
		\end{tabular}
	\end{center}
\end{table}

It follows from Table~\ref{t:mism} that $\PD(\Tfma)$ depends on the true parameter value $\theta_r$ but not on the assumed value $\theta$. Very small differences in different rows occur
only because of statistical errors of Monte Carlo simulations. This is expected since the
structure of the FMA rule does not depend on the assumed parameter value $\theta$. See \eqref{TFMA} in Remark~\ref{Rem:Monotone} where $S_t=Y_t$ in our example.
Thus, the detection probability $\PD(\Tfma)$ is only a function of $\theta_r$, but not of $\theta$. On the contrary, the probability of detection $\PD(T_{\varrho})$ varies as a function of the assumed signal intensity
$\theta$ because the structure of the modified CUSUM rule depends on the assumed value of $\theta$. Hence, the detection probability $\PD(T_{\varrho})$ is a function of both $\theta_r$ and $\theta$.

We can therefore conclude that in the sense of sensitivity with respect to the mismatch between true and assumed parameter values the FMA rule has an advantage over the modified CUSUM rule.

\subsection{Comparison of CUSUM and Modified CUSUM Rules }

The results of the comparison of the optimal modified CUSUM rule $T_{\varrho}$ against the conventional CUSUM rule $\Tcs$ are shown in
Tables \ref{t:B_gauss_inf_cusum} and \ref{t:B_gauss_inf_cusum_2}. It is seen that the operating characteristics of these rules are almost the same.
The optimal one only slightly outperforms the conventional one in all tested cases.

  We iterate that in contrast to the FMA rule structures of both rules depend on the assumed parameter value $\theta$.

\begin{table}[!ht!]	
	\renewcommand{\arraystretch}{1}
\begin{center}
	\caption{Detection characteristics of the modified CUSUM and conventional CUSUM rules ($N \sim \mrm{Geom}(\varrho),\LPFA_m = 0.001, \theta = 2.0$) \label{t:B_gauss_inf_cusum}}
	\begin{tabular}{ |c|c|c|c| }
		\hline
		\multicolumn{4}{|c|}{$m = 20$ } \\
		\hline
		$\varrho $ & 0.2 & 0.1 & 0.05 \\
		\hline
		$\PD(T_{\varrho})$ & 0.3747 & 0.6179 & 0.7855 \\
		\hline
		$\PD(\Tcs)$ & 0.3702 & 0.6121 & 0.7848 \\
		\hline
		\hline
		\multicolumn{4}{|c|}{$m = 60$ } \\
		\hline
		$\varrho $ & 0.2 & 0.1 & 0.05 \\
		\hline
		$\PD(T_{\varrho})$ & 0.3340 & 0.5787 & 0.7695 \\
		\hline
		$\PD(\Tcs)$ & 0.3286 & 0.5776 & 0.7618 \\
		\hline
		\hline
		\multicolumn{4}{|c|}{$m = 100$ } \\
		\hline
		$\varrho $ & 0.2 & 0.1 & 0.05 \\
		\hline
		$\PD(T_{\varrho})$ & 0.3216 & 0.5624 & 0.7598 \\
		\hline
		$\PD(\Tcs)$ & 0.3192 & 0.5613  & 0.7520\\
		\hline
	\end{tabular}
\end{center}
\end{table}

\begin{table}[!ht!]
\renewcommand{\arraystretch}{1}
\begin{center}
	\caption{Detection characteristics of the modified CUSUM and conventional CUSUM rules ($N \sim \mrm{Geom}(\varrho),\LPFA_m = 0.001, \theta = 1.2$) \label{t:B_gauss_inf_cusum_2}}
	\begin{tabular}{ |c|c|c|c| }
		\hline
		\multicolumn{4}{|c|}{$m = 20$ } \\
		\hline
			$\varrho $ & 0.2 & 0.1 & 0.05 \\
		\hline
		$\PD(T_{\varrho})$ & 0.1294  & 0.3392 & 0.5892  \\
		\hline
		$\PD(\Tcs)$ & 0.1271 & 0.3385 & 0.5791 \\
		\hline
		\hline
		\multicolumn{4}{|c|}{$m = 60$ } \\
		\hline
			$\varrho $ & 0.2 & 0.1 & 0.05 \\
		\hline
		$\PD(T_{\varrho})$ & 0.0927  & 0.2922  & 0.5377 \\
		\hline
		$\PD(\Tcs)$ & 0.0875 & 0.2851  & 0.5270 \\
		\hline
		\hline
		\multicolumn{4}{|c|}{$m = 100$ } \\
		\hline
			$\varrho $ & 0.2 & 0.1 & 0.05 \\
		\hline
		$\PD(T_{\varrho})$ & 0.0890  & 0.2728  & 0.5163 \\
		\hline
		$\PD(\Tcs)$ & 0.0839 & 0.2656 & 0.5066\\
		\hline
	\end{tabular}
\end{center}
\end{table}

\section{Application to Detection of Streaks of Space Objects}\label{sec:Appl}

\begin{figure}[!h!] \centering
	\includegraphics[height=0.2\textheight, width = 0.45\textwidth]{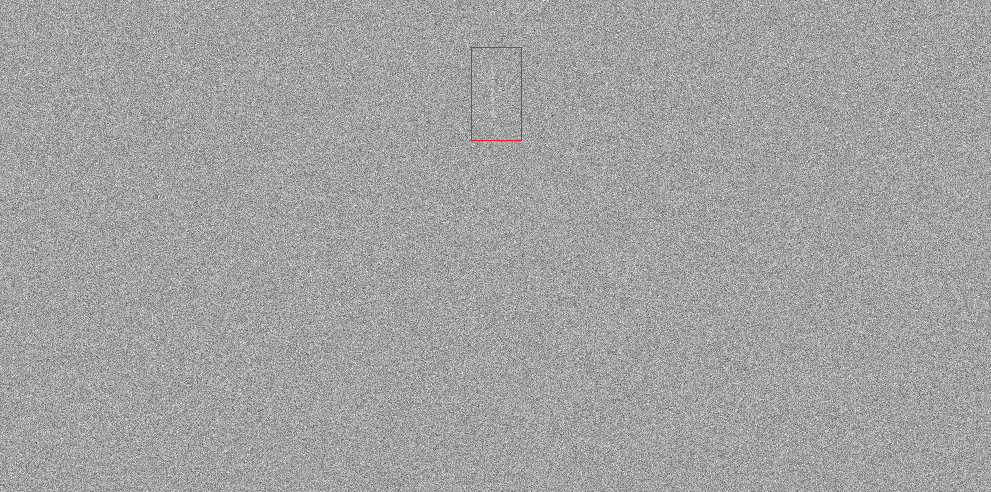}
	\caption{Digital frame with a low-contrast streak (SNR $\approx 1$). Red rectangle marks the streak position.}
	\label{fig:fit_track}
\end{figure}

Extracting streaks of faint space objects with unknown orbits from digital frames, captured with ground-based telescopes, is an important problem for Space Informatics.
A variety of methods is employed using intra-frame as well as inter-frame data processing (see, e.g., \cite{RadimSara_2019},\cite{Hickson2018},\cite{Kolessa_2013},\cite{Kol_Ravd_Pruglo},\cite{Defence2007}).
In this section, the FMA procedure, described in Section \ref{sec:ADP}, is applied to the detection of a satellite streak with a low signal-to-noise ratio (SNR) in the digital frame (see Fig. \ref{fig:fit_track}).

We use a two-stage approach. As our problem is similar to the changepoint detection problem (since the distribution of observations changes abruptly when the streak starts and ends), at the first stage, we use the FMA procedure \eqref{FMA}, which will be re-written for the 2-D space case. Since FMA detects streaks with a delay (i.e., we can only localize streaks), at the next stage, we use maximum likelihood estimation to calculate the streak position more accurately.

When operating with real frames situation is aggravated by the presence of stars and background that produce strong discrete clutter. In this case, special preprocessing for clutter removal have to be implemented (see, e.g., \cite{Tartakovsky&Brown-IEEEAES08}). For our problem, we consider that the input frame is free from clutter and contains only one streak and Gaussian noise after preprocessing, independent from pixel to pixel.

To be specific, let the satellite have a linear and uniform motion in the frame. The satellite streak is given by the vector $\Xb=(x_0,y_0,x_1,y_1)$, where $(x_0,y_0)$ corresponds to the start point and $(x_1, y_1)$ corresponds to the endpoint.

\begin{figure}[!h!] \centering
	\includegraphics[height=0.25\textheight, width= 0.5\textwidth]{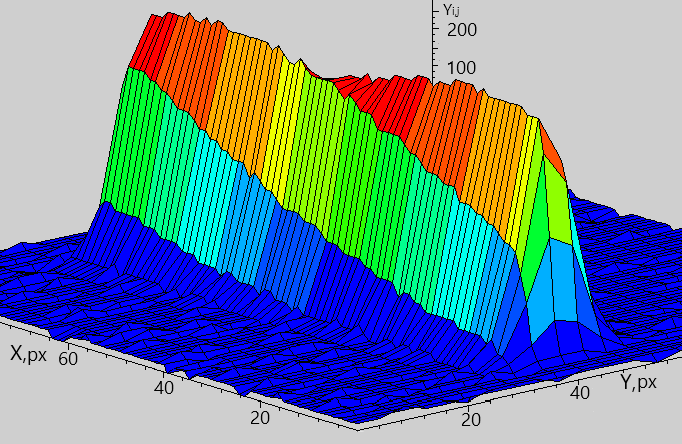}
	\caption{Model profile of the streak with the length of 80 px.} 
	\label{fig:signal}
\end{figure}

Hence, consider the following model of the observation $Y_{i,j}$ in pixel $(i,j)$ of the 2-D frame \cite{Kolessa}:
\begin{equation}\label{eq:signal_main}
Y_{i,j}= AS_{i,j}(\Xb) + \epsilon_{i,j},
\end{equation}
where $A$ is a signal intensity from the object,
$\{S_{i,j}(\Xb)\}$ are values of the model profile of the streak that are calculated beforehand assuming point spread function (PSF) is Gaussian with a certain effective width (see Fig.~\ref{fig:signal}); and
$\epsilon_{i,j}$ is Gaussian noise after preprocessing with zero mean and known (estimated empirically) local variance $\sigma^2$. Thus, the observation $Y_{i,j}$ has normal (``pre-change'') distribution $g(Y_{i,j}) = \Nc(0,\sigma^2)$
when the streak does not ``cover" pixel $(i,j)$ and normal (``post-change'') distribution $f(Y_{i,j}) = \Nc(AS_{i,j}(\Xb),\sigma^2)$ when the streak ``covers" the pixel $(i,j)$.

The problem is to detect the streak with the maximal probability of detection as well as to find the most accurate estimate $\hat{\Xb} = (\hat{x_0},\hat{y_0},\hat{x_1},\hat{y_1})$ of the streak position or
to make a decision that streak does not exist in the frame.


\begin{remark}
It is worth noting that the problem of extraction of objects' steaks can be  solved non-sequentially as a fixed sample size joint hypothesis testing and estimation problem, using, e.g., 
the generalized likelihood ratio hypothesis (GLR) test. However, the GLR test  requires testing too many hypotheses since neither direction nor position of streaks is known. As a result, the GLR test is usually 
very time-consuming. The computational complexity of the proposed changepoint detection (with further estimation) algorithm is quite low. This is an mportant advantage over the GLR test.
\end{remark}

\subsection{Stage 1: Detection and Localization of the Streak}

\begin{figure}[!h!] \centering
	\includegraphics[ height=0.2\textheight, width = 0.45\textwidth]{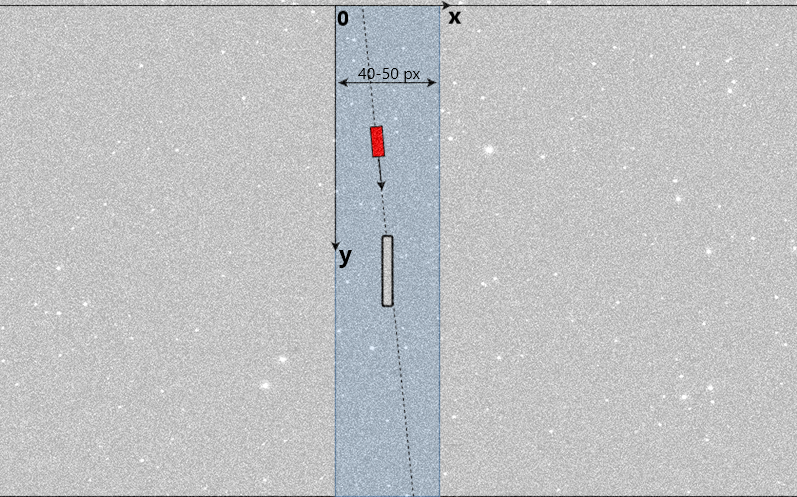}
	\caption{Search area is blue. Dotted line shows one of the possible directions. White rectangle -- streak, red rectangle -- sliding window.}
	\label{fig:search_area}
\end{figure}

In what follows, we restrict our attention to 
 intra-frame detection of faint streaks of only subequatorial satellites with unknown orbits on frames taken with telescopes mounted at the equator. Thus, a signal from each satellite (with unknown intensity, start and end points) is located almost vertically in a small area at the center of the frame shown in Fig. \ref{fig:search_area}.

Let $\Omega_S$ denote streak search area. We select a step of $0.5$ px in the upper and lower borders of the search area $\Omega_S$ to define a set of directions inside $\Omega_S$. Let $d$ denote a certain direction. Let $M_d(k),k \ge 1$  stand for a 2-D sliding rectangular window which contains certain pixel numbers $(i,j)$ at each step $k$ in the direction $d$.
Window $M_d(k)$ has a fixed length of $N_d$ pixels and a fixed width of $K_d$ pixels (the choice of the parameters depends on the expected SNR and PSF effective width).

\begin{figure}[!h!] \centering
	\includegraphics[ height=0.2\textheight, width = 0.25\textwidth]{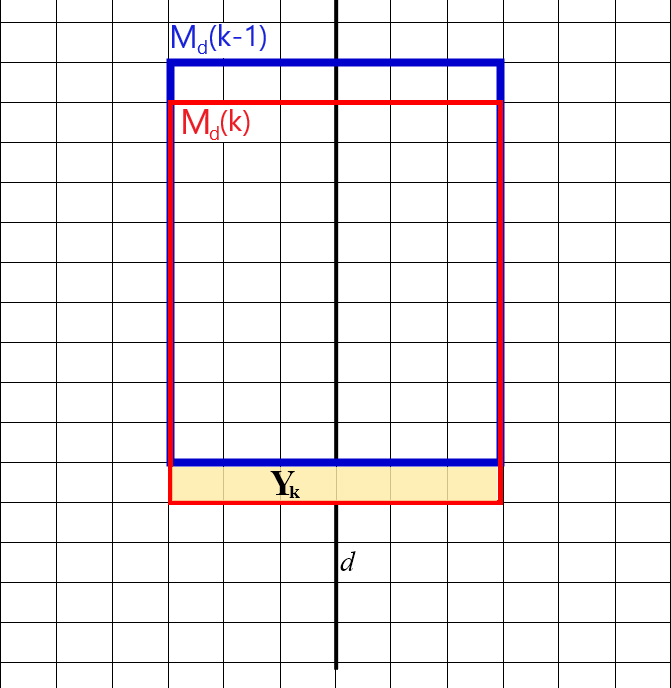}
	\caption{Observations $\{\Yb_k\}_{k \ge 2}$ contain signals $Y_{i,j}$ from new pixel numbers $(i,j)$ which appear in the window $M_d(k)$ when sliding from the previous $(k-1)$-th step.}
	\label{fig:observ}
\end{figure}

If we fix the certain direction $d$, then for this direction we solve the sequential detection problem.
Thus, observations $\{Y_n\}_{n \ge 1}$ from Section~\ref{sec:Model} are formed in the following manner (see Fig. \ref{fig:observ}):
\[
\begin{split}
&\Yb_1 ~\text{is formed from signals}~\{Y_{i,j}\}_{(i,j) \in M_d(1)},\\
&\{\Yb_k\}_{k \ge 2} ~\text{are formed from signals} ~\{Y_{i,j}\}_{(i,j) \in M_d(k) - M_d(k-1)} .
\end{split}
\]

For the Gaussian model of the streak profile and considering that $k$ plays the role of the time index $n$, defined in Sections~\ref{sec:Model} and \ref{sec:Rule}, the FMA procedure \eqref{FMA} is re-written as
\[
\widehat{\Tfma}(\tilde{h})=\inf\left\{k \ge 1:R_{M_d(k)}(\Yb) \ge \tilde{h} \right\},
\]
\[
R_{M_d(k)}(\Yb) = \sum_{(i,j)\in M_d(k)}Y_{i,j}S_{i,j}(\Xb_1),
\]
where $\{S_{i,j}(\Xb_1)\}$ are values of the Gaussian model profile calculated beforehand. Profile location is given by the vector $\Xb_1 = (x^0_1,y^0_1,x^1_1,y^1_1)$; points $(x^0_1,y^0_1)$
and $(x^1_1,y^1_1)$ are located in the direction $d$;
$R_{M_d(k)}(\Yb)$ is the changepoint detection statistic.

Experimentally we have chosen the most suitable parameters of the window $M_d(k)$ in our case: $N_d = 15$ px and $K_d = 8$ px, while SNR $\approx$ 1 and the length of the streak $N$ is unknown and 
cannot be less than $20$ px.

Hence, when sliding the 2-D window in various directions inside $\Omega_S$ and then choosing the longest sequence of $R_{M_d(k)}(\Yb)$ values above the threshold we determine the approximate position of the streak with a typical accuracy of 5-10 pixels. Therefore, we can determine an area (localization area), which with a high probability contains the streak.

\begin{figure}[!h!] \centering
	\includegraphics[height=0.22\textheight, width = 0.35\textwidth]{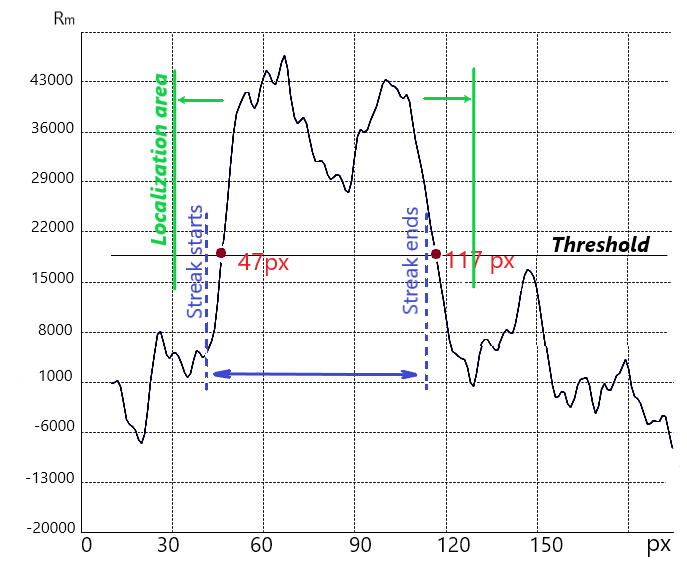}
	\caption{ The behavior of the $R_{M_d(l)}(\Yb)$ statistic along the correct direction. Real streak position is marked with blue. The streak is detected with coordinates of start and end at points 47 and 117, respectively, while the true values are 40 and 110.}
	\label{fig:plot_stat_R}
\end{figure}

Fig. \ref{fig:plot_stat_R} shows one of the possible behaviors of the $R_{M_d(k)}(\Yb)$ statistic along the right direction in the case of a very low SNR $ = 0.9$.
A more accurate estimation of the streak position is performed in the localization area at Stage~2.

\subsection{Stage 2. Accurate Estimation of the Streak Position}

After localizing the streak in a certain area (denote it by $\Pi_1$), we use maximum likelihood estimates $\hat{A}$ and $\hat{\Xb}$ of the streak parameters $\Xb,A$ calculated from
\[
(\hat{\Xb},\hat{A}) = \argmin_{\Xb,A}\sum_{(i,j)\in \Pi_1}[Y_{i,j}-AS_{i,j}(\Xb)]^2,
\]
where minimization is under constraints $(x_0,y_0) \in \Pi_1$,~~$(x_1,y_1) \in \Pi_1,~~ A > 0$ . 
It can be easily seen that
\[
\hat{\Xb}= \argmin_{\Xb}\sum_{(i,j)\in \Pi_1}\brcs{Y_{i,j} - \frac{\displaystyle\sum_{(i,j)\in \Pi_1} Y_{i,j}S_{i,j}(\Xb)}{\displaystyle\sum_{(i,j)\in \Pi_1}S_{i,j}(\Xb)^2}S_{i,j}(\Xb)}^2 .
\]
This maximum likelihood estimate $\hat{\Xb}$ of the streak position yields the final solution to our problem.

\subsection{Testing}
We tested the proposed algorithm using Monte Carlo simulation of random linear streaks with the length of 50 px superimposed on Gaussian noise.

\begin{figure}[!h!] \centering
	\includegraphics[height=0.22\textheight, width = 0.45\textwidth ]{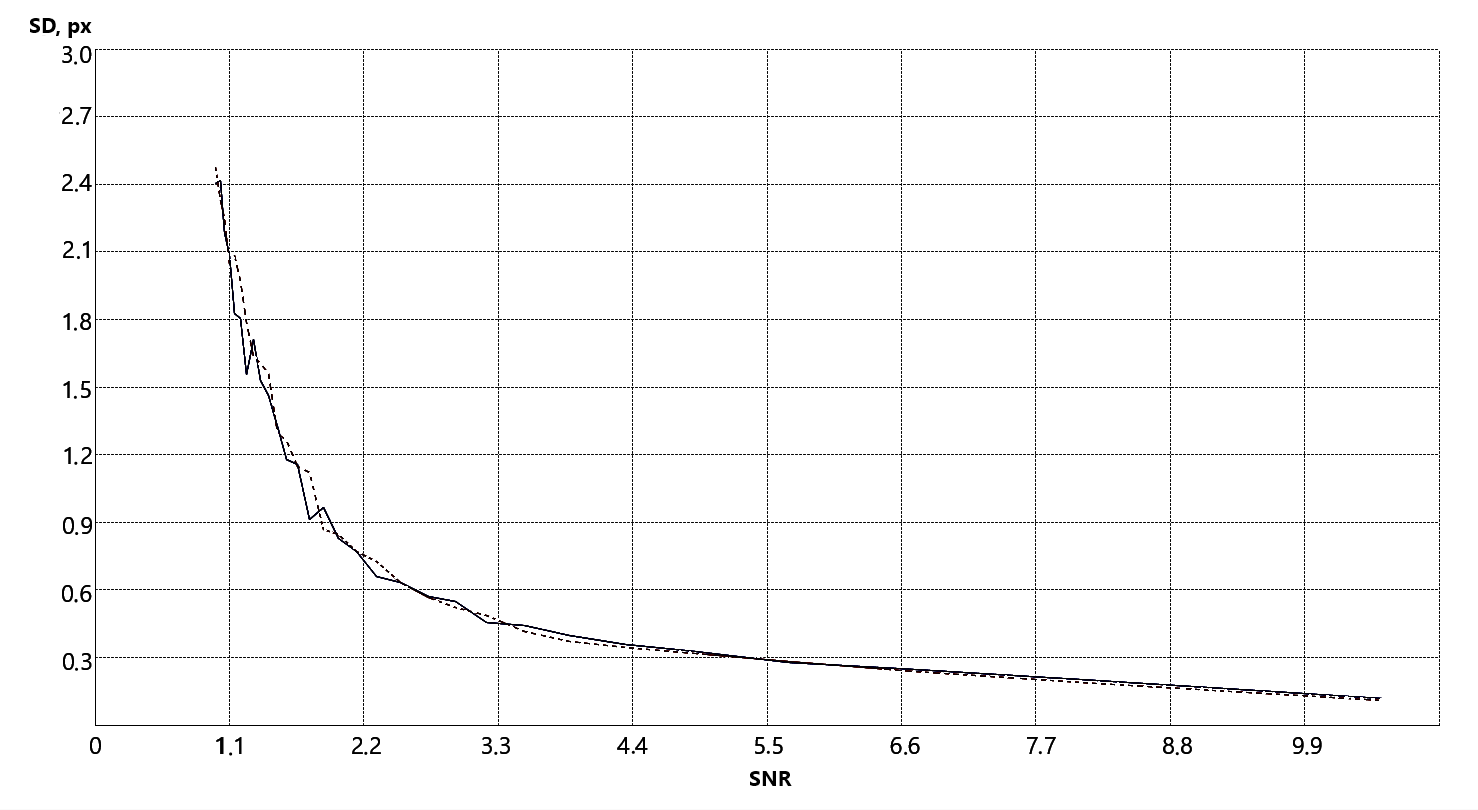}
	\caption{SD as function of SNR. Dotted line corresponds to the plot for SD of the estimate of the streak start, solid line -- to the estimate of the streak end. Number of MC trials $5\cdot10^5$.}
	\label{fig:SD(SNR)_1}
\end{figure}

The standard deviation (SD) of the estimated start and endpoints of the streak as a function of SNR was obtained for the range of SNR values from $10$ to $0.9$ (see Fig. \ref{fig:SD(SNR)_1}).
It is seen that the accuracy is high even for the SNR =1.

\section{Conclusion}\label{sec:Remarks}

We have found a strictly optimal solution to the problem of detecting signals with unknown moments of appearance and disappearance when it is required to maximize the probability of detection in the
window of random size, distributed according to geometric distribution, in the class of detection procedures with the given maximal probability of a false alarm in a prescribed
finite window. This solution is obtained using the optimal stopping theory, and the optimal procedure is a modified CUSUM. The results of Monte Carlo simulations for the Gaussian example
show that the operating characteristics of the optimal procedure are typically close to that of the more practical 
FMA procedure. As a result, we propose to use the latter procedure in an important for Space Informatics problem  -- for
intra-frame detection of faint satellite streaks with unknown orbits in optical images.

\section*{Acknowledgement}

We are grateful to George Moustakides for useful discussions.
Thanks also go to referees and Greg Sokolov whose comments improved the article.


\renewcommand{\theequation}{A.\arabic{equation}}
\setcounter{equation}{0}

\section*{Appendix}

To prove Theorem~\ref{Th:MaximinOptimalityMm} we need the following lemma, which establishes that class $\class(m,\alpha)$ is more stringent than $\class_\gamma=\{T: \Eb_\infty[T]\ge \gamma\}$ for
some appropriately selected $\gamma=\gamma(m,\alpha)$. The unconditional variant of this lemma has been considered in the paper by Lai \cite{LaiIEEE98}.

\begin{lemma}\label{Lem:LCPFAtoARL}
Let $m$ be a positive integer. If  $T\in\class(m, \alpha)$, then $T$ necessarily belongs to class $\class_\gamma$ for some $\gamma=\gamma(m,\alpha)$, in particular for
\begin{equation}\label{gammambeta}
\gamma(m,\alpha) = \frac{3}{2} + \frac{m}{\alpha} \brc{1 - \frac{3}{2} \alpha}.
\end{equation}
\end{lemma}

\begin{IEEEproof}
Let $m$ be a positive integer, $m <\gamma$. If $T\in\class_\gamma$, then there exists an $\ell$, possibly depending on $\gamma$, such that
\begin{equation} \label{prob0}
\Pb_\infty(T \le \ell+m | T > \ell)<m/\gamma .
\end{equation}
Indeed, we have
\begin{align} \label{ARLneat}
\Eb_\infty [T]  &= \sum_{\ell=0}^{\infty} \Pb_\infty\brc{T > \ell} =\sum_{i=0}^{m-1}\sum_{k=0}^{\infty} \Pb_\infty\brc{T > i+k m} \nonumber
\\
 &= \sum_{i=0}^{m-1}\sum_{k=0}^{\infty}\Pb_\infty(T > i)\Pb_\infty(T > i+k m | T > i) .
\end{align}
Suppose that
\begin{equation*}
\Pb_\infty(T > \ell+m | T> \ell)<1-m/\gamma \quad \text{for all $\ell \in\Zbb_+$}.
\end{equation*}
Then
\[
\Pb_\infty(T > i + k m | T> i)< (1-m/\gamma)^k,
\]
which can be easily derived by induction over $j=1,\dots,k$ from the equality  
\begin{align*}
&\Pb_\infty(T > i + j m | T> i)  =  \Pb_\infty\set{T>i+ (j-1) m | T>i} 
\\
& \times \Pb_\infty\set{T>i+(j-1)m +m | T > i+(j-1)m}  . 
\end{align*}
So it follows from equality \eqref{ARLneat} that
\begin{align*}
\Eb_\infty [T]  &<  \sum_{i=0}^{m-1}\Pb_\infty\brc{T > i} \sum_{k=0}^{\infty}(1-m/\gamma)^k
\\
&= (\gamma/m)\sum_{i=0}^{m-1}\Pb_\infty\brc{T > i} < \gamma,
\end{align*}
which contradicts the assumption $T\in\class_\gamma$, and therefore, proves  \eqref{prob0}.

Next, we prove that, for a given $0<\alpha < 1$, the constraint
\begin{equation} \label{LCPFAconstr}
\sup_{\ell \in \Zbb_+} \Pb_\infty(T \le \ell+m | T > \ell) \le \alpha \quad \text{for some $m \ge 1$}
\end{equation}
is stronger than the average run length constraint $\Eb_\infty [T] \ge \gamma$ ($\gamma \ge 1$), i.e., if $T\in \class(m,\alpha)$, then this implies that $T\in\class_\gamma$ for some
$\gamma=\gamma(\alpha,m)$. If the inequality \eqref{LCPFAconstr} holds, then 
\[
\Pb_\infty(T > i + k m | T> i) \ge (1-\alpha)^k  \quad \text{for all} ~ i \ge 0, 
\]
and using \eqref{ARLneat} we obtain
\begin{align*}
\Eb_\infty [T] & = \sum_{i=0}^{m-1} \Pb_\infty(T > i) \sum_{k=0}^{\infty}\Pb_\infty(T > i+k m | T > i) 
\\
& \ge  \sum_{i=0}^{m-1}\Pb_\infty\brc{T > i} \sum_{k=0}^{\infty}(1- \alpha)^k =
\\
&=  \frac{1}{\alpha} \sum_{i=0}^{m-1}\Pb_\infty\brc{T > i}.
\end{align*}
Now, obviously, for any $j=1,2,\dots$ and $i\in ((j-1)m, jm]$
\begin{align*}
\Pb_\infty(T \le i) & \le \sum_{k=1}^j \Pb_\infty\set{(k-1) m< T \le km} 
\\
& \le j \alpha < \brc{1+ \frac{i}{m} }\alpha
\end{align*}
since for any $T\in \class(m,\alpha)$
$$
\Pb_\infty\set{(k-1) m< T \le km} \le \Pb_\infty\set{T > (k-1) m} \alpha \le \alpha .
$$
It follows that 
\begin{align*}
& \sum_{i=0}^{m-1}\Pb_\infty\brc{T > i}  = 1+ \sum_{i=1}^{m-1} \brcs{1-\Pb_\infty\brc{T \le i}} 
 \\
 & \ge 1+ (m-1)(1-\alpha) -\frac{\alpha}{m} \sum_{i=1}^{m-1} i = \frac{3}{2} \alpha + \brc{1- \frac{3}{2}\alpha} m ,
\end{align*}
and therefore, 
\[
\Eb_\infty [T] \ge  \frac{3}{2} + \frac{m}{\alpha} \brc{1 - \frac{3}{2} \alpha}.
\]
This implies \eqref{gammambeta}.
\end{IEEEproof}

\begin{IEEEproof}[Proof of Theorem~\ref{Th:MaximinOptimalityMm}]
We have
\begin{align*}
&\sum_{i=1}^\infty \pi_i \Pb_\nu(\nu < T \le \nu+i |T>\nu, \Yb_1^\nu, M=i)
\\
&= \Eb_\nu \brcs{\sum_{i=1}^\infty \varrho(1-\varrho)^{i-1} \Ind{0< T-\nu \le i}| T>\nu, \Yb_1^\nu, M=i}
\\
& =  \Eb_\nu \brcs{\sum_{i=T-\nu}^\infty \varrho(1-\varrho)^{i-1}|T>\nu, \Yb_1^\nu}
\\
&= (1-\varrho)^{-1} \Eb_\nu\brcs{(1-\varrho)^{T-\nu}|T>\nu, \Yb_1^\nu}.
\end{align*}
As a result, the original maximin problem \eqref{MaximinPDOpt1Mrandom_1} reduces to
\begin{equation}\label{OP1}
\inf_{T\in \class(m,\alpha)} \inf_{\nu\in \Zbb_+} \esssup \, \frac{1}{\varrho} \Eb_\nu\brcs{1-(1-\varrho)^{T-\nu} | T>\nu, \Yb_1^\nu} .
\end{equation}

Consider the optimization problem
\[
\inf_{T\in\class_\gamma} \inf_{\nu\in \Zbb_+} \esssup \, \frac{1}{\varrho} \Eb_\nu\brcs{1-(1-\varrho)^{T-\nu} | T>\nu, \Yb_1^\nu}
\]
in class $\class_\gamma$.  By Theorem~2.1 in~\cite{poor-as98}  the optimal solution for this problem is the stopping time $T_{\varrho}(B)$ with threshold $B=B_\gamma$ such that
$\Eb_\infty [T_{\varrho}(B)]=\gamma$. But by Lemma~\ref{Lem:LCPFAtoARL} class $\class_\gamma \subset\class(m,\alpha)$ for some $\gamma=\gamma(m,\alpha)$, so that the
stopping time $T_{\varrho}(B)$ with threshold $B=B(m,\alpha)$ that satisfies equality \eqref{Bmalpha} and $\Eb_\infty[T_{\varrho}(B(m,\alpha))]=\gamma(m,\alpha)$ solves the optimization problem \eqref{OP1}, which implies that
\begin{align*}
&\sup_{T\in \class(m,\alpha)}\inf_{\nu\in \Zbb_+} \essinf\, \bar{\Pb}_\nu(T \le \nu+M | \Yb_1^\nu, T >\nu)
\\
& = \sup_{T\in \class_{\gamma(m,\alpha)}} \inf_{\nu\in \Zbb_+} \essinf\, \bar{\Pb}_\nu(T \le \nu+M | \Yb_1^\nu, T >\nu)
\\
&  =\inf_{\nu\in\Zbb_+} \essinf\, \bar{\Pb}_\nu(T_{\varrho} \le \nu+M | \Yb_1^\nu, T_{\varrho} >\nu),
\end{align*}
where $T_\varrho=T_{\varrho}(B(m,\alpha))$, and the proof is complete.
\end{IEEEproof}

%




\begin{IEEEbiography}[{\includegraphics[width=1.2in,height=1.45in,keepaspectratio]{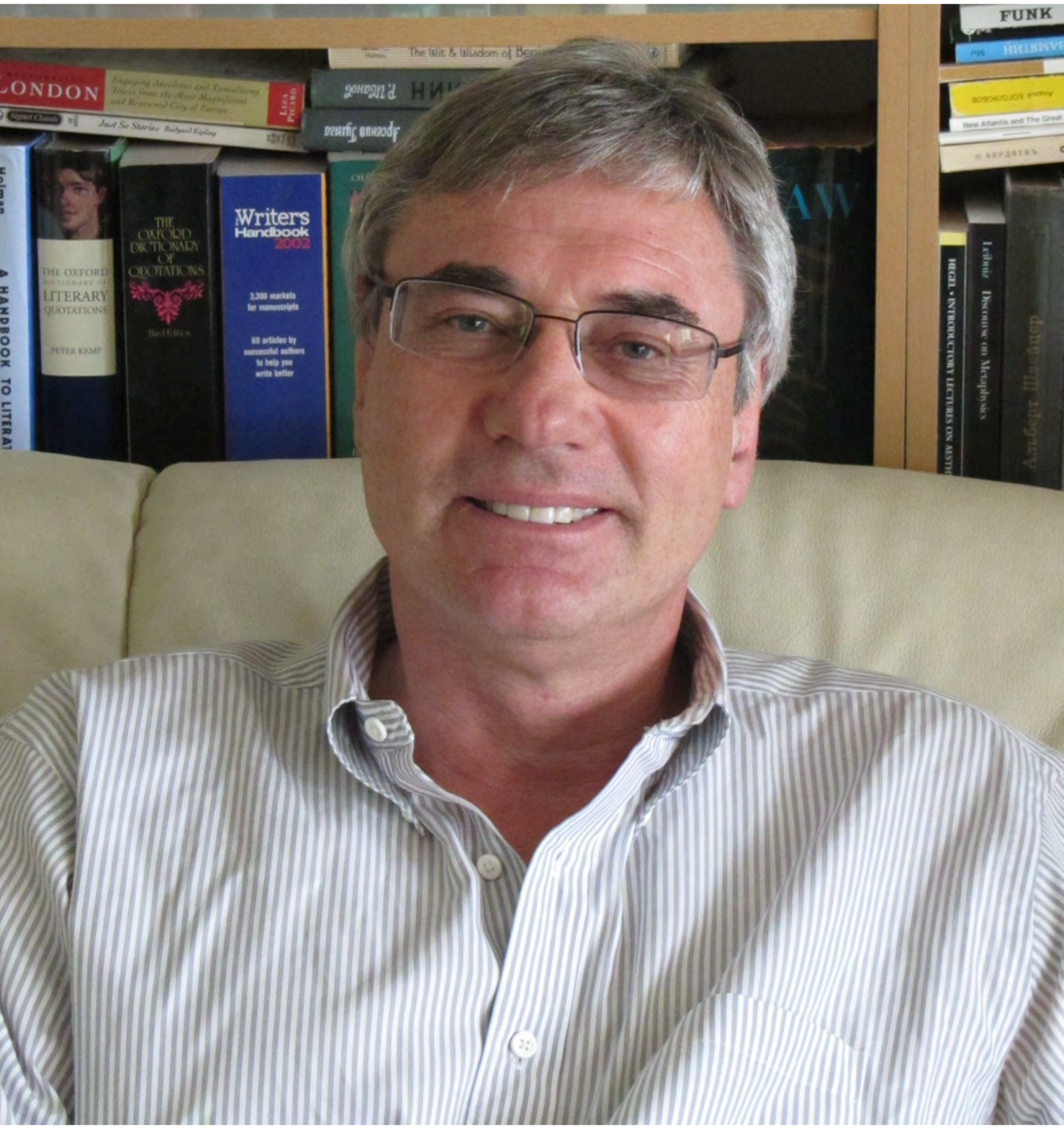}}]{Alexander G.\ Tartakovsky} (M'01-SM'02),  M.S., Ph.D., D.Sc., is an award-winning statistician, 
Head of the Space Informatics Laboratory at the Moscow Institute of Physics and Technology 
(``PhysTech''), and President of AGT StatConsult, Los Angeles, CA. From 2013 to 2015, he was a Professor of Statistics at the University of Connecticut, Storrs. Previously, for almost two decades, 
he was a Professor in the Department of Mathematics and the Associate Director of the Center for Applied Mathematical Sciences at the University of Southern California (USC).

Dr.\ Tartakovsky is the author of three books, several book chapters, and over 100 papers across a range of subjects, including theoretical and applied statistics; applied probability; sequential analysis; and changepoint detection. 
His research has many applications, including in statistical image and signal processing, video tracking, detection and tracking of targets in radar and infrared search and track systems, near-Earth space informatics, 
information integration/fusion, intrusion detection and network security, rapid detection of epidemics, and detection and tracking of malicious activity.

Dr.\ Tartakovsky earned an M.S. in Electrical Engineering from the Moscow Aviation Institute in 1978 and a Ph.D. in Statistics and Information Theory from PhysTech in 1981. 
He also earned an advanced Doctor of Science (D.Sc.) degree from PhysTech in 1990.

From 1981 to 1992, he was first a Senior Research Scientist and then Department Head at the Moscow Institute of Radio Technology and a Professor at PhysTech, 
where he worked on the application of statistical methods to optimization and modeling of information systems.

From 1993 to 1996, Dr.\ Tartakovsky was a professor at the University of California, Los Angeles (UCLA), first in the Department of Electrical Engineering and then in the Department of Mathematics.

Dr.\ Tartakovsky has received numerous awards for his work, including the Abraham Wald Prize in Sequential Analysis and several Best Young Scientist awards from the Russian Academy of Sciences. 
He is also a Fellow of the Institute of Mathematical Statistics (IMS) and a senior member of the Institute of Electrical and Electronics Engineers (IEEE).
\end{IEEEbiography}

\begin{IEEEbiography}[{\includegraphics[width=1.2in,height=1.4in,keepaspectratio]{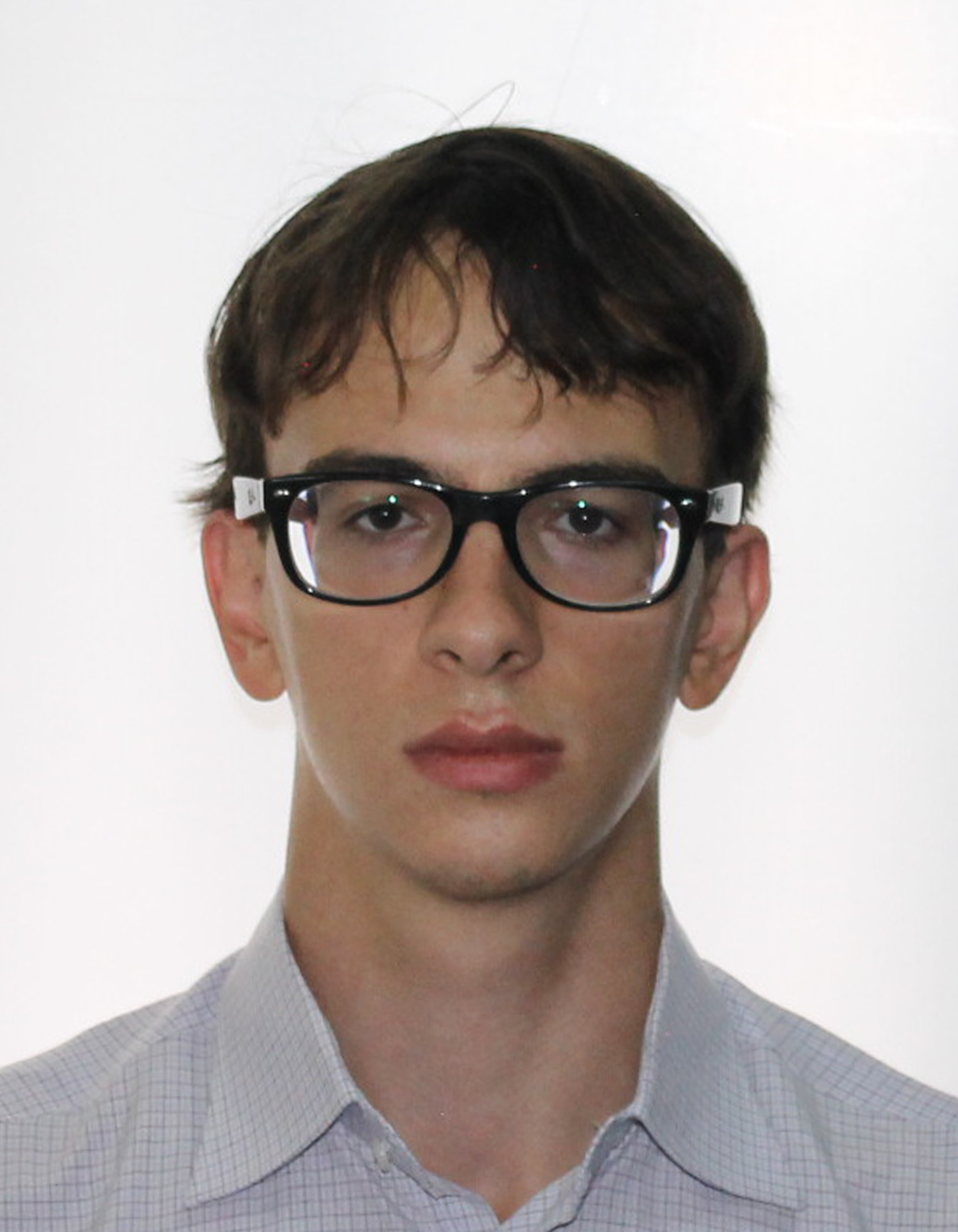}}]{Nikita R.\ Berenkov} is a postgraduate student at the Moscow Institute of Physics and
Technology (MIPT), Russia and engineer in the Space informatics Laboratory at MIPT.
\end{IEEEbiography}

\begin{IEEEbiography}[{\includegraphics[width=1.2in,height=1.4in,keepaspectratio]{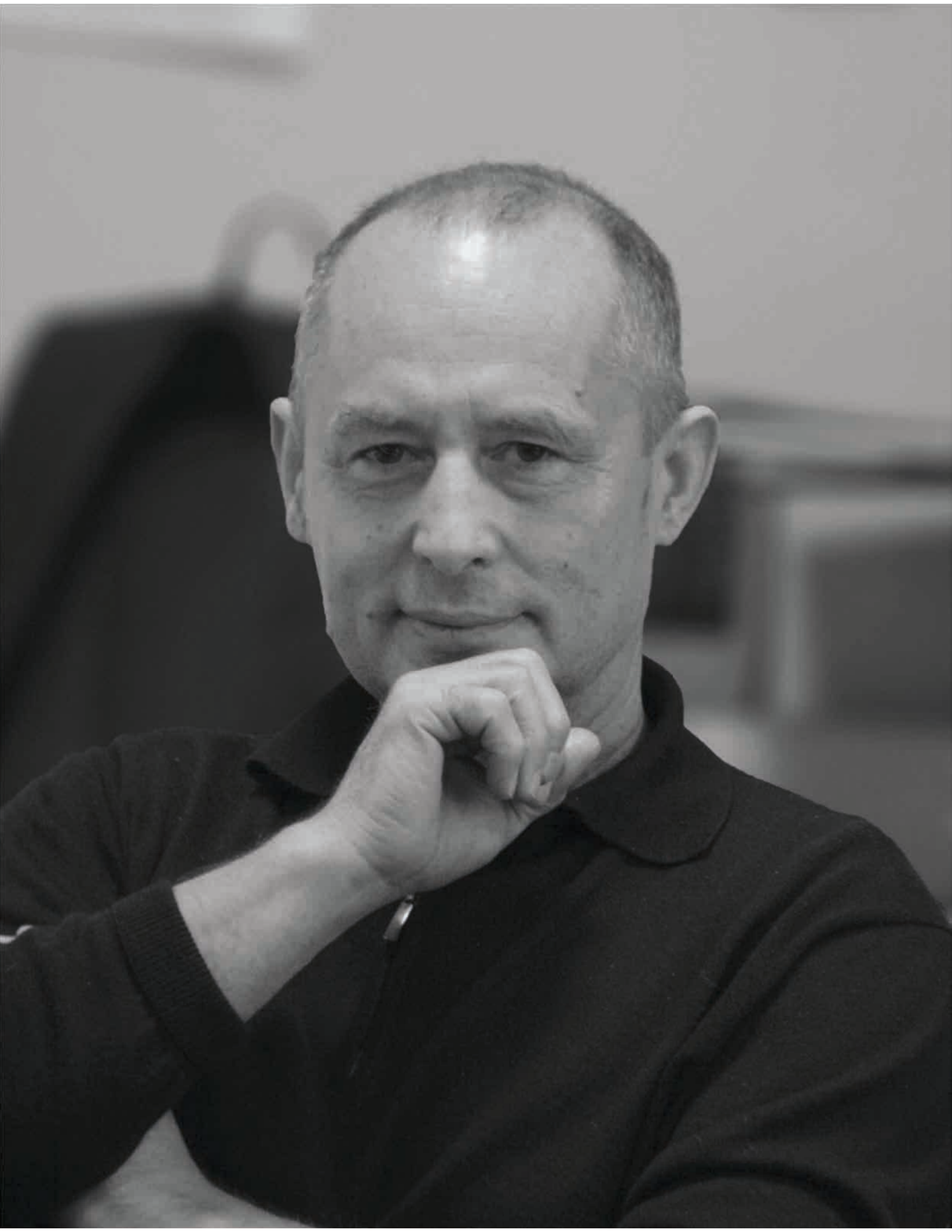}}]{Alexey E.\ Kolessa}
research interests include nonlinear filtering, object detection and tracking, video surveillance, image and signal processing, near-Earth space informatics.
He has over than 50 publications. Dr. Kolessa is the Principal Scientist in the Space informatics Laboratory at the Moscow Institute of Physics and Technology (MIPT), Russia as well
as the deputy chair of the Information Systems Department, MIPT.
\end{IEEEbiography}

\begin{IEEEbiography}[{\includegraphics[width=1.2in,height=1.4in,keepaspectratio]{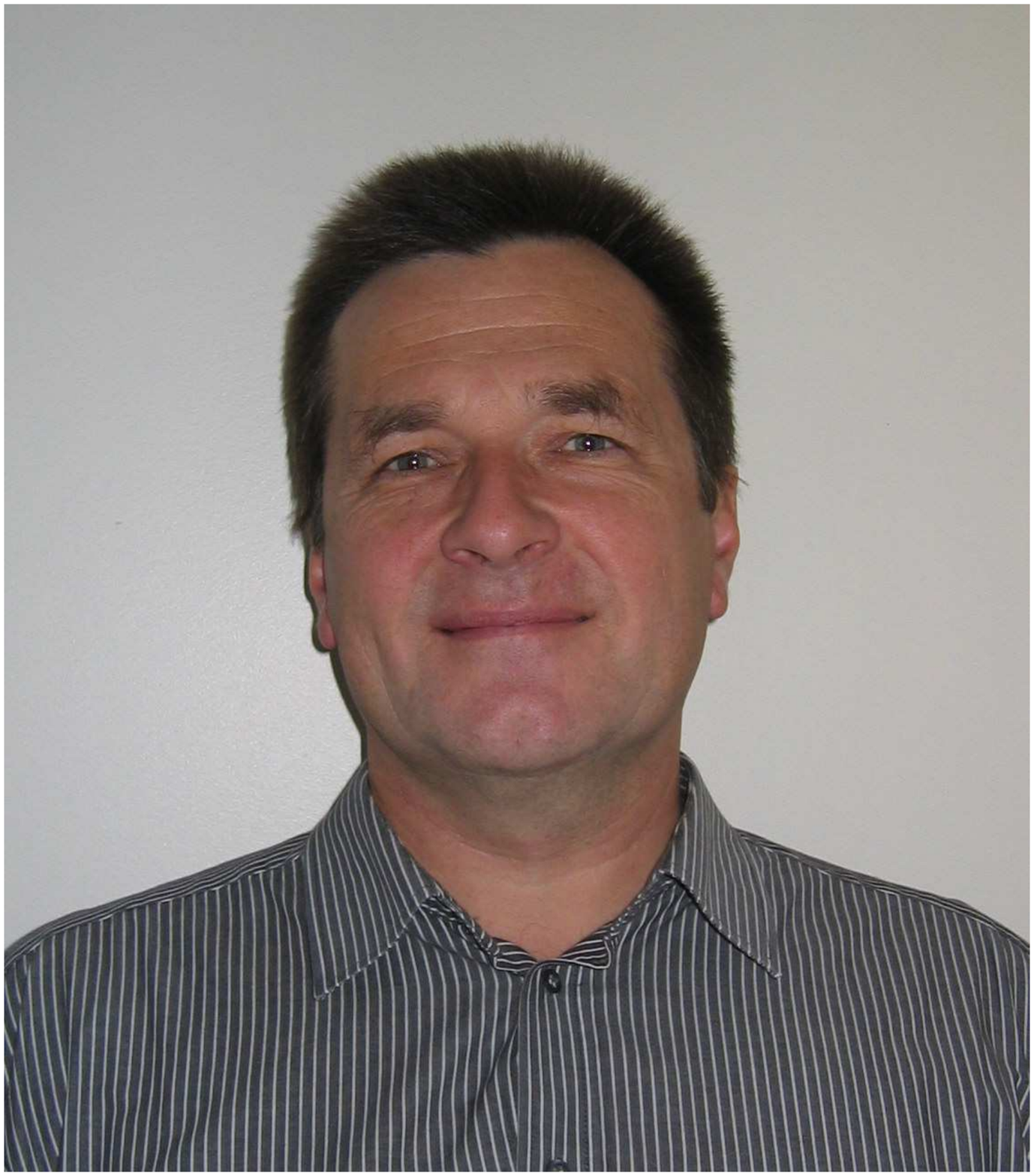}}]{Igor V. Nikiforov} is Emeritus Professor at the
Universit\'e de Technologie de Troyes, France. He is currently active or has demonstrated activity in the past in the following areas: sequential change detection and isolation (multiple hypotheses case); statistical signal detection with nuisance parameters; statistical fault detection and isolation; sequential analysis; multivariable control theory and statistical process control. The results of his research have found implementations in the areas of navigation; hidden information/channel detection; safety-critical system monitoring; network anomaly detection; digital parametric tomography; seismic signal processing and tsunami warning systems. He is the author/co-author of five books, 18 chapters, 56 journal papers and 183 conference papers/presentations. He received a {\em 2017 Abraham Wald Award in Sequential Analysis.}

Igor Nikiforov received his M.S. degree in automatic control from Moscow Institute of Physics and Technology in 1974, and the Ph.D. in automatic control from the Institute of Control Sciences (USSR Academy of Science), Moscow, in 1981. He spent 1974 -- 1992 as a research engineer, junior scientist, senior scientist, head of research group at the Institute of Control Sciences. He spent as a senior scientist 1992 -- 1994 at the IRISA/INRIA, France and 1994 -- 1995 at the University of Lille, France. He joined, as professor, the University of Technology of Troyes (UTT) in 1995. He was the head of laboratory LM2S (UTT) during 2000 -- 2007 and the head of the Institute of Computer Sciences and Information of Troyes (UTT), FRE CNRS 2732, during  2004 -- 2005.
\end{IEEEbiography}

\vfill

\end{document}